\documentclass[12pt]{article}
\usepackage{amsmath, amsthm, amssymb}
\usepackage{xspace}
\usepackage{cite}
\usepackage{fancyhdr}
\newcommand{\comment}[1]{}

\newcommand{\trivialexamples}[1]{#1}
\theoremstyle{plain}
\newtheorem{theorem}{Theorem}[section]
\newtheorem{lemma}[theorem]{Lemma}

\newtheorem*{theorem*}{Theorem}
\newtheorem*{lemma*}{Lemma}

\theoremstyle{definition}
\newtheorem{definition}[theorem]{Definition}
\newtheorem*{definition*}{Definition}
\newtheorem{notation}[theorem]{Notation}
\newtheorem*{notation*}{Notation}
\newtheorem{example}[theorem]{Example}
\newtheorem*{example*}{Example}
\theoremstyle{remark}
\newtheorem{remark}[theorem]{Remark}
\newtheorem*{remark*}{Remark}

\let\oldenumerate\enumerate
\renewcommand{\enumerate}{
   \oldenumerate
\renewcommand{\labelenumi}{\textit{(\arabic{enumi})}}
}

\DeclareMathOperator{\dom}{dom}
\DeclareMathOperator{\len}{len}

\title{An extension of the \ef game for \mbox{first
order} logics augmented with \lind quantifiers}
\date{}
\author{Simi Haber\\
        Saharon Shelah}

\newcommand{\edge}{\sim}
\newcommand{\ef}{Ehrenfeucht-Fra\"{\i}sse\xspace}
\newcommand{\lind}{Lindstr\"{o}m\xspace}
\newcommand{\modelClassType}[1]{\ensuremath{\mathbf{#1}}}
\newcommand{\taugra}{\ensuremath{\tau_\text{\sc{Gra}}}\xspace}

\begin{document}
\maketitle
\thispagestyle{fancy}

\begin{abstract}
We propose an extension of the \ef game able to deal with logics augmented with \lind quantifiers.
We describe three different games with varying balance between simplicity and ease of use.
\end{abstract}

\emph{Dedicated to Yuri Gurevich on the occasion of his 75th birthday}

\section{Introduction}
The \ef game \cite{A:fraisse50, phd:fraisse53, A:fraisse54, A:ehrenfeucht61} is an important tool in contemporary model theory, allowing to determine whether two structures are elementary equivalent up to some quantifier depth. It is one of the few model theoretic machineries that survive the transition from general model theory to the finite realm. \comment{It also makes the field more accessible to combinatorialists, who are stereotypically fond of analyzing game strategies, and may be less familiar with the ins and outs of logic.}

There are quite a few known extensions of the \ef game and in the following we mention a few (this is not a comprehensive list). In \cite{B:immerman1998} Immerman describes how to adapt the \ef game in order to deal with finite variable logic. Infinitary logic has a precise characterization by a similar game \cite{A:barwise1977, A:immerman1982}. An extension for fixpoint logic and stratified fixpoint logic was provided by Bosse \cite{IP:bosse1993}.

\lind quantifiers were first introduced and studied by \lind in the sixties \cite{A:lindstrom1966a, A:lindstrom1966b, A:lindstrom1966c, A:lindstrom1969} and may be seen as precursors to his theorem.

The aim of this paper is to present several related extensions of the \ef game adapted to logics augmented with \lind quantifiers.

\section{The Game}
\begin{notation} \begin{enumerate}
\item Let $\tau$ denote a vocabulary. We assume $\tau$ has no function symbols, but that is purely for the sake of clearer presentation. $\tau$ may have constant symbols.
\item First order logic will be denoted by $\mathcal{L_{FO}}$. Along this paper we will look on extensions of first order logic, therefore the logic under discussion will change according to our needs. We shall denote the logic currently under discussion by $\mathcal{L}$, and we will explicitly redefine $\mathcal{L}$ whenever needed.
\item Given a vocabulary $\tau$, we use $\mathcal{L}(\tau)$ to denote the \emph{language} with logic $\mathcal{L}$ and vocabulary $\tau$. We will use this notation only when clarity demands, so in fact we may abuse notation and use $\mathcal{L}$ also for the language under discussion.
\item For even further transparency, all the examples in this manuscript (in particular, all cases of pairs of models to be proved equivalent) will be dealing with simple\footnote{An undirected graph with no loops and no double edges is called \emph{a simple graph}.} graphs. Hence (only in examples) we further assume that $\tau$ is the vocabulary of graphs denoted henceforth by \taugra. Explicitly, $\taugra = \{ \edge \}$ where $\edge$ is a binary, anti-reflexive and symmetric relation. For the \lind quantifiers given in examples, we may use vocabularies other than \taugra.

\item Let $\modelClassType{A}_1, \modelClassType{A}_2, \dotsc$ be classes of models\footnote{I don't see a problem with having infinitely many of these.}, each closed under isomorphism. The models in $\modelClassType{A}_i$ are all $\tau_i$-structures in some relational vocabulary
$\displaystyle\tau_i = \left\{P_{i,1}^{a_{i,1}},\dotsc,P_{i,t_i}^{a_{i,t_i}}
\right\}$.
\item For simplicity, we will assume that each $\tau_i$ has an additional relation, $P_{i,0}^1$. This will serve for the formula defining the universe of the model. Formally, all our models will have their domain the entire universe, and the first relation will be a subset defining the domain de facto.
\item We set $a_{i,0}=1$ for every $i$.
\item Each $\modelClassType{A}_i$ corresponds to a \lind quantifier $Q_i$ binding $a_i = \sum_{j=0}^{t_i} a_{i,j}$ variables.
\end{enumerate}
\end{notation}

%\footnote{Should I include $x_0$ in the count? right now it is out so the last sentence is missing by 1.}.%
\trivialexamples{
\begin{example}\begin{enumerate}
\item $\modelClassType{A}_1$ may be the class of commutative groups, in which case $\tau_1$ is consisted of a constant symbol and a ternary relation encoding the group operation.
\item Another example may be finite Hamiltonian graphs, in which case the vocabulary is the vocabulary of graphs and the class $\modelClassType{A}$ will be the set of all finite Hamiltonian graphs (over, say, $[n]=\{1,\dots,n\}$ for any $n\in\mathbb{N}$).
\end{enumerate}\end{example}
}%

\begin{notation}
Giver a vector $\bar{x}$, we denote its length by $\len(\bar{x})$.
\end{notation}
\begin{definition} \label{def: lind quant satisfaiability}
We define the quantifier $Q_i$ corresponding to $\modelClassType{A}_i$ as follows: Let $G$ be a $\tau$-structure with domain $V$. For any index $i$ and formulae
$\varphi_0(x_0, \bar{y}), \varphi_1(\bar{x}_1, \bar{y}), \dotsc,
\varphi_{t_i}(\bar{x}_{t_i}, \bar{y})$ such that $\len(\bar{x}_j) = a_{i,j}$, the satisfiability of $Q_i\, x_0, \bar{x}_1, \dotsc, \bar{x}_{t_i} (
\varphi_0(x_0, \bar{b}), \varphi_1(\bar{x}_1, \bar{b}), \dotsc,
\varphi_{t_i}(\bar{x}_{t_i}, \bar{b}))$ is given by
\begin{align*}
G \models Q_i\, &x_0, \bar{x}_1, \dotsc, \bar{x}_{t_i} (
\varphi_0(x_0, \bar{b}), \varphi_1(\bar{x}_1, \bar{b}), \dotsc,
\varphi_{t_i}(\bar{x}_{t_i}, \bar{b}) ) \Longleftrightarrow \\ %
( &\{x_0 \in V \mid G \models \varphi_0(x_0, \bar{b})\},
\{\bar{x}_1\in V^{a_{i,1}} \mid G \models \varphi_1(\bar{x}_1,
\bar{b})\}, \dotsc, \\ %
& \{\bar{x}_{t_i} \in V^{a_{i,t_i}} \mid G
\models \varphi_{t_i}(\bar{x}_{t_i}, \bar{b})\}) \in \modelClassType{A}_i,
\end{align*}
where $\bar{b}$ are parameters.
\end{definition}

\begin{remark}
Definition \ref{def: lind quant satisfaiability} requires $\varphi_0$ to have exactly one free variable, $x_0$ (excluding $\bar{y}$, saved for parameters). However there is not real reason to to avoid sets of vectors of any length from serving as the domain of the model defined in the quantifier. We will not discuss this here, but the generalization of the proposed games to this case are straightforward.
\end{remark}

\begin{definition}
\begin{enumerate}
\item Let $\tau$ be a vocabulary and $\mathcal{L}=\mathcal{L}(\tau)$ be a language. Given two $\tau$-structures $G_1, G_2$ (not necessarily with distinct universe sets) and two equal length sequences of elements $\bar{x}_1 \in G_1$, $\bar{x}_2 \in G_2$, we say that $(G_1,\bar{x}_1)$ and $(G_2,\bar{x}_2)$ are \emph{$k$-equivalent with respect to $\mathcal{L}$} if for any formula $\varphi(\bar{x}) \in \mathcal{L}$ of quantifier depth at most $k$ one has 
\[
G_1 \models \varphi(\bar{x}_1) \Longleftrightarrow G_2 \models \varphi(\bar{x}_2).
\]
\item When considering only one model, that is, when we take $G=G_1=G_2$, we refer to the equivalence classes of this relation in the domain of $G$ simply by \emph{the $(a,k,G)$-equivalence classes} (or just equivalence classes when the context is clear enough).
\end{enumerate}
\end{definition}
Notice that unions of $(a,k,G)$-equivalence classes are exactly the definable sets of $a$-tuples of elements in $\dom(G)$ using $\mathcal{L}$-formulas of quantifier depth at most $k$.
\trivialexamples{%
\begin{example}
Let $\mathcal{L}$ be the first order language of graphs, $\mathcal{L}=\mathcal{L_{FO}}(\taugra)$, and let $G = (V,E)$ be a graph. If $G$ is simple then the $(1,0,G)$-equivalence classes are $V$ and $\varnothing$. If $|V|>1$ then the $(1,1,G)$-equivalence classes are\footnote{The atomic sentences appearing in $\varphi(x)$ are $x=y$ and $x\edge y$.} the set of isolated vertices in $G$, the set of vertices adjacent to all other vertices and the set of vertices having at least one neighbor and one non-neighbor (some of which may be empty of course).
\end{example}
\begin{notation}
We denote the logic obtained by augmenting the first order logic with the quantifiers $Q_1, Q_2, \dotsc$ by $\mathcal{L} = \mathcal{L}[Q_1, Q_2, \dotsc]$.
\end{notation}
\begin{example}
Consider $\mathcal{L} = \mathcal{L}[Q_\text{\sc{Ham}}](\taugra)$, where $Q_\text{\sc{Ham}}$ stands for the ``Hamiltonicity quantifier'' (corresponding to the class of graphs containing a Hamiltonian cycle --- a cycle visiting each vertex precisely once). Let $G$ be a graph. Then the set of all vertices $x$ for which all of the graphs\footnote{Here $N(x) = \{y \in V \mid x\edge y\}$ is the neighborhood of $x$ in $G$, $G[U]$ where $U\subseteq V$ is the graph induced on $U$ and $\overline{G}$ is the compliment of $G$. } $G[N_G(x)], \overline{G}[N_G(x)], G[N_{\overline{G}}(x)], \overline{G}[N_{\overline{G}}(x)]$
are Hamiltonian is an example of a $(1,1,G)$-equivalence class with respect to $\mathcal{L}[Q_\text{\sc{Ham}}]$. The set of vertices with degree exactly two is a union of $(1,1,G)$-equivalence classes, as can be seen by\footnote{$\varphi$ expresses: ``the complete graph $K_{d(x)}$ is Hamiltonian'' which is true when $d(x)>2$ and false when $d(x)=2$ (we may treat $K_0$ and $K_1$ separately, if needed).} 
\[
\varphi(x) = Q_\text{\sc{Ham}} x_0,x_1,x_2(\, x_0 \edge x, x_1 \neq x_2).
\]
\end{example}
}%

\newcommand{\efl}{\ensuremath{\mathrm{EFL}}\xspace}
\newcommand{\dup}{\ensuremath{\mathrm{ISO}}\xspace}
\newcommand{\spo}{\ensuremath{\mathrm{AIS}}\xspace}

\newcommand{\efla}{\ensuremath{\efl_1}\xspace}
\newcommand{\eflb}{\ensuremath{\efl_2}\xspace}
\newcommand{\eflc}{\ensuremath{\efl_3}\xspace}

\subsection{Description of the first game}
Before describing the game, we need the following definition:
\begin{definition} \label{def: copy of M in G}
Let $\tau$ be a vocabulary, $\mathcal{L}$ a language over that vocabulary (not necessarily first order) and $G$ a model of $\tau$. Additionally, let $M=(S',R_1',\dots,R_{t}')$ be a model of another vocabulary $\tau'$. A \emph{copy of $M$ in $G$} is a tuple $(S,R_1,\dots,R_t)$ such that
\newcounter{contdEnumi}
\begin{enumerate}
\item $S$ is a subset of $\dom(G)$ with the same cardinality as $\dom(M)=S'$ (where $\dom(G)$ is the universe or underlying set of $G$).
\item $R_1,\dots,R_t$ are relations over $S$, such that each $R_j$ has the same arity as $R_j'$.
\item $(S,R_1,\dots,R_t)$ is isomorphic to $(S',R_1',\dots,R_{t}')$.
\setcounter{contdEnumi}{\value{enumi}}
\end{enumerate}

If in addition the following holds
\begin{enumerate} \setcounter{enumi}{\value{contdEnumi}}
\item \label{def part: induced condition} $S$ is union of $(1,k,G)$-equivalence classes, and each relation $R_j$ of arity $a_j$ is a union of $(a_j,k,G)$-equivalence classes;
\end{enumerate}
we say that a copy of $M$ in $G$ is \emph{$k$-induced by $\mathcal{L}$}. When $k$ and / or $\mathcal{L}$ can be clearly determined by the context, we may omit mentioning one of them, or both.
\end{definition}

We may now define the first game.
\begin{definition}
Let $G_1$ and $G_2$ be two models with domains $V_1$ and $V_2$ respectively. Let $k\geq 0$ an integer and $\bar{c}_\ell = (c^1_\ell, \dotsc, c^r_\ell) \in V_\ell^r$ two finite sequences. We define the game\footnote{We will describe a few variants, hence the subscript.} $\efla[G_1, G_2, \bar{c}_1, \bar{c}_2\, ;\, k]$. There are two players, named \dup and \spo. The game board is the models $G_1$ and $G_2$ plus the sequences $\bar{c}_\ell$ and there are $k$ rounds. Each round is divided into two parts, and each part consists of two sub-rounds. The game is defined recursively. If $k = 0$, then if the mapping $c_1^i \rightarrow c_2^i$ is an isomorphism, then $\dup$ wins, otherwise $\spo$ wins.

When $k>0$ then first \spo plays. He picks one of the models $G_1$
or $G_2$ (denoted henceforth by $G_\ell$) and a quantifier $Q_i$ (or
the existential quantifier\footnote{In this case, $\modelClassType{A}_\exists = P(V)\setminus\{\varnothing\}$, so \spo may choose any non-empty subset $S_\ell$ of $V_\ell$.}). Next \spo picks a model $M \in \modelClassType{A}_i$,
and embeds it into $G_\ell$ in a manner that preserve $(k-1,G_\ell)$-equivalence classes. That is, \spo picks a tuple $(S_\ell,R_{\ell,1},\dots,R_{\ell,t_i})$ that is a copy of $M$ in $G$ which is $(k-1)$-induced by $\mathcal{L}$ enriched with $r$ constants having values $\bar{x}_\ell$.
\comment{%
of $V_\ell$ of the same size as the universe of $M$ and provides relations $R_{\ell,1},\dots,R_{\ell,t_i}$ over $S_\ell$ such that:
\begin{enumerate}
\item $(S_\ell, R_{\ell,1},\dots,R_{\ell,t_i})$ is isomorphic to $M$; and
\item each of the set $S_\ell$ and the relations $R_{\ell,j}$ is a union of $(a_{i,j},k-1)$-equivalence classes with respect to $\mathcal{L}$ enriched with $r$ constants having values $\bar{c}_\ell$
\end{enumerate}
We call this process ``picking an induced copy of $M$ in $G_\ell$''.
}%
If \spo can not find such an embedding, he loses\footnote{We will consider only logics stronger than first-order, hence the existential quantifier is always assumed to be at \spo' disposal and he will never lose in this manner.}. Implicitly \spo claims that \dup can not find a matching induced copy of a model from $\modelClassType{A}_i$.

Second, \dup responds by choosing a model $M'$ from $\modelClassType{A}_i$ ($M'$ may not necessarily be the same as $M$), and then picking an induced copy of $M'$ in $G_{3-\ell}$ which we naturally denote by $(S_{3-\ell}, R_{3-\ell,1}, \dots, R_{3-\ell,t_i})$. She is implicitly claiming that her choices match the picks of \spo, that is, each $R_{3-\ell,j}$ (or $S_{3-\ell,j}$) is a union of $(a_{i,j},k-1,G_{3-\ell})$-equivalence classes defined by the same formulas as the formulas defining the $(a_{i,j},k-1,G_{\ell})$-equivalence classes of which $R_{\ell,j}$ is made. If \dup can not complete this part she loses. This ends the first part of the round.

In the second part of the round \spo chooses $m\in\{1,2\}$ and $0 \leq j \leq t_i$. He then picks $(c_m^{r+1}, \dotsc, c_m^{r+a_{i,j}}) \in R_{m,j}$ (implicitly challenging \dup to do the same). Finally \dup picks $(c_{3-m}^{r+1}, \dotsc, c_{3-m}^{r+a_{i,j}}) \in R_{3-m,j}$ and they move on to play
\begin{align*}
\efla[G_1, G_2, & (c_1^1, \dotsc, c_1^r, c_1^{r+1}, \dotsc,
c_1^{r+a_{i,j},}), \\
& (c_2^1, \dotsc, c_2^r, c_2^{r+1}, \dotsc, c_2^{r+a_{i,j},}); k-1] .
\end{align*}
This ends the second part and the round. Since $k$ goes down every round, the game ends when $k=0$, as described above.
\end{definition}

Given the description above, the following should be self-evident:
\begin{lemma} \label{lem: efl1 is useful}
Let $\mathcal{L}=\mathcal{L}[Q_1,Q_2,\dots](\tau)$ be a language over some vocabulary $\tau$ where $Q_1,Q_2,\dots$ are \lind quantifiers and let $G_1, G_2$ be two $\tau$-structures. Then, $\dup$ has a winning strategy for $\efla[G_1, G_2, \varnothing, \varnothing\, ; \,
k]$ if and only if for any sentence $\varphi \in \mathcal{L}$ of
quantifier depth at most $k$
\[ G_1 \models \varphi \Longleftrightarrow G_2 \models \varphi . \]
\end{lemma}
\comment{
\begin{proof}
When $k=0$, the quantifier free formulas are simply conjunctions, disjunctions and negation of the atomic formulas, the latter containing only the constants as variables. Hence we can assert isomorphism of the structure induced on the constants, but noting more. This is identical with the winning condition of \dup.

Assume $k>0$ and the validity of the Lemma until $k-1$.
\end{proof}
}
\subsection{A game where definability is not forced}
While the claim of Lemma \ref{lem: efl1 is useful} is satisfying, it may be hard to put into use since it takes finding unions of $(a,k-1,G)$-equivalence classes for granted, being a rule of the game. This might hinder strategy development and we would like to describe another game with looser rules, denoted \eflb.

In this version the players are not bound to choosing unions of $(a,k-1,G)$-equivalence classes when picking a copy of the chosen model (hence we call their action ``picking a copy of $M$ in $G_\ell$'', omitting the ``induced'' part). That is, we omit requirement \ref{def part: induced condition} in Definition \ref{def: copy of M in G}. It falls to the other player to check that indeed every relation is a union of the relevant equivalence classes. A general round now goes as follows:

\spo picks a graph $G_\ell \in \{G_1,G_2\}$ and a quantifier $Q_i$ (or, as before, the existential quantifier). Next \spo picks a model $M \in \modelClassType{A}_i$ and picks a copy of $M$ in $G_\ell$. His implicit claim now includes the claim that each of the relations he chose is a union of $(a_{i,j},k-1,G_\ell)$-equivalence
classes with respect to $\mathcal{L}$ enriched with $r$ constants
having values $\bar{c}_\ell$.

$\dup$ can respond in two different ways --- she can ``accept the
challenge'' (as she did in \efla), or attack the second part of the claim of \spo. That is, she can do one of the following:
\begin{enumerate}
\item\label{enumi: dup accepts integrity}%
Accept. In this case she chooses $M'\in \modelClassType{A}_i$ and picks a copy of $M'$ in $G_{3-\ell}$.
Implicitly she is claiming that her choices matches the choices of
\spo. \trivialexamples{That is, the set of vertices $S_{3-\ell}$ and each of the relations defined on it are a union of the $(a_{i,j},k-1,G_{3-\ell})$-equivalence
classes corresponding\trivialexamples{\footnote{We say that $E_1$, an $(a,k,G_1)$-equivalence class of $a$-tuples in $G_1$ corresponds to $E_2$ --- a set of $a$-tuples in $G_2$ if for any $\bar{x}_1 \in E_1$ and $\bar{x}_2 \in E_2$ one has \[ G_1 \models \varphi(\bar{x}_1) \Leftrightarrow G_2 \models \varphi(\bar{x}_2)\] for any $\varphi \in \mathcal{L}$ of quantifier depth at most $k$.}} to the ones that $\spo$ picked.} This ends the first part of the round.

\spo may continue in a two different ways.
\begin{enumerate}
\item\label{enumii: spo rejects integrity}%
Reject the fact that $S_{3-\ell}$ or one of the relations picked by
$\dup$ is a union of equivalence classes. In order to settle this, we recursively use \eflb:\newline
Again, we let $a=a_{i,j}$ be the arity of
the allegedly invalid relation $R_{3-\ell,j}$. $\spo$ picks two $a$-tuples,
$(c^{r+1}, \dotsc, c^{r+a}) \in R_{3-\ell,j}$ and $(c'^{r+1}, \dotsc,
c'^{r+a}) \in V_{3-\ell}^a\setminus R_{3-\ell,j}$, and they move on to
play
\begin{align*}
\eflb[G_{3-\ell}, G_{3-\ell}, & (c_{3-\ell}^1, \dotsc, c_{3-\ell}^r,
c^{r+1}, \dotsc, c^{r+a}), \\ %
& (c_{3-\ell}^1, \dotsc, c_{3-\ell}^r, c'^{r+1}, \dotsc, c'^{r+a});
k-1] .
\end{align*}
with exchanged roles (since this time \spo claims the two tuples are actually $(a,k-1,G_{3-\ell})$-equivalent).

\item\label{enumii: spo accepts integrity}%
Reject the fact that \dup's choice matches his choice (as he did in \efla). In this
case he picks a relation $P_j \in \tau_i$ and an $a_{i,j}$-tuple of
elements from $S_\ell$ (or one element if he challenges her choice
of $S_{3-\ell}$). Denote his choice by $(c_\ell^{r+1}, \dotsc,
c_\ell^{r+a_{i,j}}) \in S_\ell$. $\dup$ responds by picking another
$a$-tuple $(c_{3-\ell}^{r+1}, \dotsc, c_{3-\ell}^{r+a_{i,j},}) \in
S_{3-\ell}$, and they move on to play
\begin{align*}
\eflb[G_1, G_2, & (c_1^1, \dotsc, c_1^r, c_1^{r+1}, \dotsc,
c_1^{r+a_{i,j},}), \\
& (c_2^1, \dotsc, c_2^r, c_2^{r+1}, \dotsc, c_2^{r+a_{i,j},}); k-1] .
\end{align*}
\end{enumerate}

\item\label{enumi: dup rejects integrity}%
Reject. In this case $\dup$ wants to prove that $S_\ell$ or one of
the relations picked by $\spo$ is not a union of equivalence
classes. We continue similarly to case 1.(b):\newline
Let $a=a_{i,j}$ be the arity
of the allegedly relation, $R_{\ell,j}$, splitting an equivalence class. $\dup$ picks two $a$-tuples,
$(c^{r+1}, \dotsc, c^{r+a}) \in R_{\ell,j}$ and $(c'^{r+1}, \dotsc,
c'^{r+a}) \in V_\ell^a\setminus R_{\ell,j}$, and they move to play
\begin{align*}
\eflb[G_{\ell}, G_{\ell}, & (c_{\ell}^1, \dotsc, c_{\ell}^r,
c^{r+1}, \dotsc, c^{r+a}), \\ %
& (c_{\ell}^1, \dotsc, c_{\ell}^r, c'^{r+1}, \dotsc, c'^{r+a});
k-1] .
\end{align*}
this time keeping their original roles.
\end{enumerate}

For any two models $G_1$ and $G_2$, constants $\bar{c}_1, \bar{c}_2$ and $k \in \mathbb{N}$, whoever has a winning strategy for $\efla[G_1,G_2,\bar{c}_1, \bar{c}_2;k]$ has a winning strategy for $\eflb[G_1,G_2,\bar{c}_1, \bar{c}_2;k]$. Hence the parallel of Lemma \ref{lem: efl1 is useful} is true for \eflb as well.

While we got the benefit of in-game validation of the equivalence classes integrity claims, \eflb is not easy to analyze in applications because the game-board and players role change over time. We amend this in the last suggested version of the game.

\subsection{A game with fixed game-board and roles}
The last version, denoted \eflc, forks from \eflb in two places.

\begin{definition}
We define \eflc like \eflb except that:
\begin{enumerate}
\item First, assume the game reaches step \ref{enumi: dup rejects integrity}., where \dup wants to prove that \spo's chose a relation $R_{\ell,j}$ splitting an equivalence relation. In this case the first part of the round ends immediately and the second part goes as follows:\newline
\dup chooses two $a_{i,j}$-tuples, $\bar{c}_{\ell,1}$ from $R_{\ell,j}$ and $\bar{c}_{\ell,2}$ from the complement of $R_{\ell,j}$. She then pick another $a_{i,j}$-tuples from $G_{3-\ell}$, denoted $\bar{c}_{3-\ell}$. Spoiler than picks one of $\bar{c}_{\ell,1}$ or $\bar{c}_{\ell,2}$ and they move on to play \eflc with $\bar{c}_{3-\ell}$ concatenated to the constants of $G_{3-\ell}$ and \spo's choice concatenated to the constants of $G_{\ell}$, and $k-1$ moves. They keep their roles and the game-board is still $G_1$ and $G_2$.

If \dup can not find a matching tuple in $G_{3-\ell}$, she can not disprove the integrity claim of \spo, but it does not matter as $G_1$ and $G_2$ are not $k$-equivalent and she is bound to lose anyway.

Notice that in this case the first part of the round had only \spo playing, and in the second part \dup played first.

\item The second (and last) change from \eflb happens when the game is in step \ref{enumii: spo rejects integrity}. In this case \spo wants to prove that \dup's choice of at least one relation $R_{3-ell,j}$ is splitting an equivalence relation. In this case \spo picks a tuple $\bar{c}_{3-\ell}$ (from the suspicious equivalence class) in $G_{3-\ell}$ that is not in $R_{3-ell,j}$ and challenges \dup to find a matching tuple $\bar{c}_{\ell}$ in $G_{3-\ell}$ that is not in $R_{ell,j}$. They move on to play \eflc with these choices and $k-1$ moves. Again both roles and game-board remain as was. Notice that the game flow in this case is actually the same as the game flow in \ref{enumii: spo accepts integrity}.
\end{enumerate}
\end{definition}

As before, it is easy to convince oneself that the claim of Lemma \ref{lem: efl1 is useful} is still valid. We repeat it here:

\begin{lemma} \label{lem: efl3 is useful}
$\dup$ has a winning strategy for $\eflc[G_1, F_2, \varnothing, \varnothing\, ; \,
k]$ if and only if for any sentence $\varphi \in \mathcal{L}$ of
quantifier depth at most $k$
\[ G_1 \models \varphi \Longleftrightarrow G_2 \models \varphi . \]
\end{lemma}

\comment{
\section{Applications?}
}

\section{Summary}
We have presented three equivalent variants of the celebrated \ef game adapted to deal with logics extended by \lind quantifiers. We believe \eflc may be easier to analyse than direct quantifier elimination and it is out hope that it will find applications.

\comment{
\appendix
\section{map}
I'm writing this for a long time now, so here is a map of letters I use:
\begin{tabular}{cll}
  $\tau$ & base vocabulary, in our case $\{\edge\}$ &  \\
  $\edge$ & graph adjacency &  \\
  $i$ & index for \lind quantifiers \\
  $\tau_i$ & vocabulary for \lind quantifiers & \\
  $t_i$ & number of relation in $\tau_i$ & \\
  $j$ & index for relations in $\tau_i$ & \\
  $P_{i,j}$ & relation $j$ in $\tau_i$ & \\
  $a_{i,j}$ & arity of $P_{i,j}$ & \\
  $\modelClassType{A}_i$ & classes of models for \lind quantifiers & \\
  $Q_i$ & \lind quantifiers & \\
  $a_i$ & number of variables $Q_i$ binds & \\
  $G$   & a graph \\
  $V$   & vertex set \\
  $E$   & edge set \\
  $\bar{b}$ & parameter vector \\
  $x, \bar{x}$ & vertex, vertices vector \\
  $r$ & number of constant vertices in the beginning of a game or round \\
  $G_1, G_2$ & two graphs being (most of) the game board \\
  $\ell$ & index for the part of the game-board\\
  $\bar{c}_\ell$ & a vector of elements in $G_\ell$.\\
  $a$ & length of vector of vertices. related to arity. \\
  $k$ & quantifier depth or number of turns in the game \\
  $\mathcal{L}$ & a formal language \\
  $\varphi$ & formula \\
\end{tabular}
}
\bibliographystyle{abbrv}

\end{document}